\documentclass[12pt]{smfart}
\usepackage[frenchb]{babel}
\usepackage{times,amscd,verbatim}
\usepackage[T1]{fontenc} 
\usepackage[frenchb]{babel} 

\newcommand{\Z}{\mathbf{Z}}

\newcommand{\sA}{\mathcal{A}}

\newcommand{\sM}{\mathcal{M}}

\newcommand{\IM}{\operatorname{Im}}
\newcommand{\Coker}{\operatorname{Coker}}
\newcommand{\Pic}{\operatorname{Pic}}

\newcommand{\NS}{\operatorname{NS}}
\newcommand{\LAlb}{\operatorname{LAlb}}
\newcommand{\LA}[1]{\mbox{${\rm L}_{#1}{\rm Alb}$}}

\newcommand{\Tr}{\operatorname{Tr}}

\newcommand{\et}{{\operatorname{\acute{e}t}}}

\newcommand{\rg}{\operatorname{rg}}

\newcommand{\Inj}{\lhook\joinrel\longrightarrow}

\newcommand{\iso}{\overset{\sim}{\longrightarrow}}

\newcommand{\by}[1]{\overset{#1}{\longrightarrow}}

\newtheorem{thm}{Th\'eor\`eme}
\newtheorem{cor}{Corollaire}

\newtheorem{lemme}{Lemme}
\theoremstyle{definition}

\newtheorem{rque}{Remarque}

\newcounter{spec}
{\end{list}}

\begin{document}

\title{D\'emonstration g\'eom\'etrique du th\'eor\`eme de Lang-N\'eron}
\author{Bruno Kahn}
\address{Institut de Math\'ematiques de
Jussieu\\175--179 rue du Chevaleret\\75013 Paris\\France}
\date{2 mars 2007}
\email{kahn@math.jussieu.fr}
\maketitle

\begin{abstract} On donne une d\'emonstration ``sans hauteurs" du
th\'eor\`eme de Lang-N\'eron: si $K/k$ est une extension de type fini r\'eguli\`ere
et
$A$ est une $K$-vari\'et\'e ab\'elienne, le groupe $A(K)/\Tr_{K/k} A(k)$ est de type fini, o\`u
$\Tr_{K/k} A$ d\'esigne la $K/k$-trace de $A$ au sens de Chow. La m\'ethode fournit
une expression du rang de ce groupe en fonction de ceux de certains groupes de
N\'eron-Severi.
\end{abstract}

\begin{altabstract} We give a proof without heights of the Lang-N\'eron
theorem: if $K/k$ is a regular extension of finite type and 
$A$ is an  abelian $K$-variety, the group $A(K)/\Tr_{K/k} A(k)$ is finitely
generated, where $\Tr_{K/k} A$ denotes the $K/k$-trace of $A$ in the sense of Chow.
Our method computes the rank of this group in terms of certain ranks of
N\'eron-Severi groups.
\end{altabstract}

Soient $K/k$ une extension de type fini r\'eguli\`ere et $A$ une $K$-vari\'et\'e
ab\'elienne. On se propose de donner une d\'emonstration ``sans hauteurs" du
th\'eor\`eme de Lang-N\'eron:

\begin{thm}[\protect{\cite{lang2}}] \label{t1} Le groupe $A(K)/\Tr_{K/k} A(k)$ est
de type fini, o\`u
$\Tr_{K/k} A$ d\'esigne la $K/k$-trace de $A$ au sens de Chow.
\end{thm}

(La d\'emonstration de Lang et N\'eron est expos\'ee dans le langage des sch\'emas
dans \cite[App. B]{picfini}.)

On se ram\`ene imm\'ediatement au cas o\`u $k$ est alg\'ebriquement clos. Soit $X$
un mod\`ele lisse de $K/k$ choisi de telle sorte que $A$ se prolonge en un sch\'ema
ab\'elien $p:\sA\to X$. 

\begin{lemme}\label{l1} La suite
\[0\to \Pic(X)\by{p^*}\Pic(\sA)\by{j^*}\Pic(A)\to 0\]
est exacte, o\`u $j$ est l'inclusion de $A$ dans $\sA$.
\end{lemme}

\begin{proof} On a un diagramme commutatif aux lignes exactes
\[\begin{CD}
\displaystyle\bigoplus_{x\in \sA^{(1)}}\Z@>>>\Pic(\sA)@>j^*>> \Pic(A)@>>> 0\\
@A{\wr}AA @A{p^*}AA @AAA\\
\displaystyle\bigoplus_{x\in X^{(1)}} \Z@>>> \Pic(X)@>>> 0
\end{CD}\]
d'o\`u la suite exacte d\'esir\'ee, sauf l'exactitude en $\Pic(X)$. Celle-ci
r\'esulte du fait que $p$ a une section.
\end{proof}

\begin{lemme}\label{l2} Soit $A^1(\sA)\subset \Pic(\sA)$ le sous-groupe des cycles
alg\'ebriquement \'equivalents \`a z\'ero. Alors $j^*A^1(\sA)\subset \Tr_{K/k}
\hat A(k)$, o\`u
$\hat A = \Pic^0_{A/k}$ est la vari\'et\'e ab\'elienne duale de $A$.
\end{lemme}

\begin{proof}[D\'emonstration (cf. \protect{\cite[Lect. 1, lemme 1.3]{bloch}})] Par
d\'efinition, $A^1(X)$ est engendr\'e via les correspondances alg\'ebriques par
des jacobiennes de courbes, donc le lemme r\'esulte de la d\'efinition de la
$K/k$-trace.
\end{proof}

\begin{lemme}\label{l3}  On a $j^*A^1(\sA)= \Tr_{K/k} \hat A(k)$, et une suite exacte scind\'ee
\[0\to \NS(X)\by{p^*}\NS(\sA)\by{j^*}\Pic(A)/\Tr_{K/k} \hat A(k)\to 0\]
o\`u $\NS$ d\'esigne le groupe des cycles de codimension $1$ modulo l'\'equivalence
alg\'ebrique.
\end{lemme}

\begin{proof} Le lemme du serpent appliqu\'e au diagramme commutatif de suites exactes
\[\begin{CD}
0@>>> A^1(\sA)@>>> \Pic(\sA)@>>> \NS(\sA)@>>> 0\\
&& @A{p_A^*}AA @A{p^*}AA @A{p_N^*}AA\\
0@>>> A^1(X)@>>> \Pic(X)@>>> \NS(X)@>>> 0
\end{CD}\]
fournit, via le lemme \ref{l1} et compte tenu du fait que $p$ a une section, une suite exacte
\[0\to \Coker p_A^*\to \Pic A\to \Coker p_N^*\to 0.\] 

Par le lemme \ref{l2}, $\Coker p_A^*\subset \Tr_{K/k} \hat A(k)$. D'apr\`es \cite[th.
3]{picfini}, $\Coker p_N^*$ est de type fini. Le groupe $\Tr_{K/k} \hat A(k)/\Coker p_A^*$,
divisible et de type fini, est donc nul. On obtient donc un isomorphisme $\Pic(A)/\Tr_{K/k}\hat
A(k)\iso\Coker p_N^*$, comme d\'esir\'e.
\end{proof}

Le th\'eor\`eme \ref{t1} r\'esulte du lemme \ref{l3} et de la g\'en\'eration finie de
$\NS(\sA)$, d\'ej\`a utilis\'ee dans sa d\'emonstration.\qed

On obtient aussi:

\begin{cor} $\rg (\hat A(K)/\Tr_{K/k}\hat A(k)) = \rho(\sA)-\rho(X)-\rho(\hat A)$, o\`u
$\rho(Y):=\rg \NS(Y)$ pour une vari\'et\'e lisse $Y$. (Noter que $\rho(\hat A)$ est
calcul\'e ``sur $K$".)\qed
\end{cor}

\begin{rque} Cette d\'emonstration a une interpr\'etation agr\'eable en termes de $1$-motifs,
dans le cadre d\'evelopp\'e dans \cite{bvk}. Supposons seulement $k$ parfait; soit $\sM_1$ la
cat\'egorie des $1$-motifs de Deligne sur $k$ et soit $D^b(\sM_1)$ sa cat\'egorie d\'eriv\'ee
au sens de \cite[d\'ef. 1.5.2]{bvk}. Soient $\LAlb(X)$ et $\LAlb(\sA)$ les objets de
$D^b(\sM_1)$ associ\'es \`a $X$ et $\sA$ par \cite[d\'ef. 8.1.1]{bvk}, et soit $\LAlb(\sA/X)$
la fibre du morphisme $p_*:\LAlb(\sA)\to \LAlb(X)$. D'apr\`es \cite[cor. 10.2.3]{bvk}, on a,
pour toute $k$-vari\'et\'e lisse $Y$:
\[\LA{i}(Y) =
\begin{cases}
\relax [\Z[\pi_0(Y)]\to 0]& \text{si $i = 0$}\\
\relax [0\to\sA_{Y/k}^0] & \text{si $i= 1$}\\
\relax [0\to \NS_{Y/k}^*] & \text{si $i= 2$}\\
0 & \text{sinon}
\end{cases}\]
o\`u $\pi_0(Y),\sA_{Y/k}^0$ et $\NS_{Y/k}^*$ d\'esignent respectivement l'ensemble des
composantes connexes g\'eom\'etriques, la vari\'et\'e d'Albanese et le dual de Cartier du
groupe de N\'eron-Severi (au sens ci-dessus) de $Y$; les $\LA{i}(Y)$ sont calcul\'es par
rapport \`a une $t$-structure convenable sur $D^b(\sM_1)$. Les r\'esultats pr\'ec\'edents et un
calcul facile de suite exacte donnent alors
\[\LA{i}(\sA/X) =
\begin{cases}
\relax [0\to\IM_{K/k} A] & \text{si $i= 1$}\\
\relax [0\to M^*] & \text{si $i= 2$}\\
0 & \text{sinon}
\end{cases}\]
o\`u $\IM_{K/k} A$ est la $K/k$-image de $A$ et $M(\bar k):=\Pic(A_{\bar
k})/\Tr_{K/k}
\hat A(\bar k)$ est une extension 
\[0\to \hat A(\bar k K)/\Tr_{K/k} \hat A(\bar k)\to M(\bar k)\to \NS(A_{\bar k K})\to 0\]
$\bar k$ \'etant une cl\^oture alg\'ebrique de $k$.
\end{rque}

\begin{rque}
Revenons au cas o\`u $k$ est alg\'ebriquement clos, et soit $l$ un nombre premier diff\'erent
de la caract\'eristique de $k$. En jouant avec la suite spectrale de Leray relative \`a $p$, on
obtient facilement une injection
\[f:(\hat A(K)/\Tr_{K/k} \hat A(k))\otimes \Z_l\Inj H^1_\et(X,j_* T_l(\hat A)).\]

Si $k$ est la cl\^oture alg\'ebrique d'un corps $k_0$ de type fini sur son
sous-corps premier, on a \'evidemment
\[\IM f \subset \bigcup_L H^1_\et(X,j_* T_l(\hat A))^{Gal(k/L)}\]
o\`u $L$ d\'ecrit les extensions finies de $k_0$.
\end{rque}

\end{document}